\documentclass[draft,12pt]{article}

\usepackage[latin1]{inputenc}
\usepackage{amsmath,amssymb}
\usepackage{latexsym}
\usepackage{verbatim}

\newtheorem{theorem}{Theorem}[section]
\newtheorem{corollary}[theorem]{Corollary}
\newtheorem{lemma}[theorem]{Lemma}
\newtheorem{proposition}[theorem]{Proposition}
\newtheorem{definition}[theorem]{Definition}
\newtheorem{assumption}[theorem]{Assumption}
\newtheorem{remark}[theorem]{Remark}
\numberwithin{equation}{section}

\def\proof{{\medskip\noindent {\bf Proof. }}}
\def\longproof#1{{\medskip\noindent {\bf Proof #1.}}}
\def\qed{{\hfill $\square$ \bigskip}}

\def\square{{\vcenter{\vbox{\hrule height.3pt
        \hbox{\vrule width.3pt height5pt \kern5pt
           \vrule width.3pt}
        \hrule height.3pt}}}}

 \def\sB {{\cal B}} 
 \def\sE {{\cal E}} 
  
  \def\sL {{\cal L}}

\def\ol{\overline}
\def\E{{\mathbb E}}
\def\P{{\mathbb P}}
\def\norm#1{{\Vert #1 \Vert}}
\def\del{{\partial}}
\def\lam{{\lambda}}

\def\bee{\begin{equation}}
\def\bet{\begin{theorem}}
\def\bep{\begin{proposition}}
\def\bel{\begin{lemma}}
\def\bec{\begin{corollary}}
\def\bed{\begin{definition}}
\def\eee{\end{equation}}
\def\eet{\end{theorem}}
\def\eep{\end{proposition}}
\def\eel{\end{lemma}}
\def\eec{\end{corollary}}
\def\eed{\end{definition}}

\def\R{{\mathbb R}}
\def\E{{{\mathbb E}\,}}
\def\P{{\mathbb P}}

\def\lam{{\lambda}}

\def\al{{\alpha}}
\def\grad{{\nabla}}
\def\proof{{\medskip\noindent {\bf Proof. }}}
\def\longproof#1{{\medskip\noindent {\bf Proof #1.}}}
\def\qed{{\hfill $\square$ \bigskip}}

\def\eps{\varepsilon}
\def\vp{\varphi}

\def\norm#1{\Vert #1 \Vert}

 \def\qq {\qquad}
\def\del{{\partial}}

\def\ol{\overline}

\def\ni{\noindent }

\def\bs{\bigskip}

\def\square{{\vcenter{\vbox{\hrule height.3pt
        \hbox{\vrule width.3pt height5pt \kern5pt
           \vrule width.3pt}
        \hrule height.3pt}}}}

\def\tfrac#1#2{{\textstyle {\frac{#1}{#2}}}}

\def\tlint{{- \kern-0.85em \int \kern-0.2em}}  % for textstyle
\def\dlint{{- \kern-1.05em \int \kern-0.4em}}  % for displays

 \def\sB {{\cal B}} 
 \def\sE {{\cal E}} 
  
  \def\sL {{\cal L}}

\def\nn{{\nonumber}}

\def\cbetanorm#1{{\Vert #1 \Vert}_{C^\beta}}
\def\cgammanorm#1{{\Vert #1 \Vert}_{C^\gamma}}
\def\cgammaf{{\Vert f \Vert}_{C^\gamma}}
\def\cabnorm#1{{\Vert #1 \Vert}_{C^{\alpha+\beta}}}
\def\linorm#1{{\Vert #1 \Vert}_{L^\infty}}
\def\apb{{\alpha+\beta}}

\def\canorm#1{{\Vert #1 \Vert}_{C^a}}
\def\cbnorm#1{{\Vert #1 \Vert}_{C^b}}

\begin{document}

\title{Regularity results for stable-like operators}

\author{Richard F. Bass\footnote{Research partially supported by NSF grant
DMS-0601783.}}

%\date{January 5,  2005}

\maketitle

\begin{abstract}  
\noindent For $\al\in [1,2)$ we consider operators of the form
$$\sL f(x)=\int_{\R^d} [f(x+h)-f(x)-1_{(|h|\leq 1)} \grad f(x)\cdot h]
\frac{A(x,h)}{|h|^{d+\al}}$$ and for $\al\in (0,1)$ we consider the
same operator but where the $\grad f$ term is omitted. 
We prove, under appropriate conditions on $A(x,h)$, that
the solution $u$ to $\sL u=f$ will be in $C^{\al+\beta}$ if $f\in C^\beta$.

\end{abstract}

%AAAAAAAAAA

% document here

\section{Introduction}\label{S:intro}

Many models in mathematical physics, financial mathematics,
and mathematical economics are based on  diffusions corresponding
to second order elliptic differential operators. 
In the last decade or so,
though, researchers in these areas have found that frequently real world
phenomena are better fitted if one allows jumps. 
To give a very simple example, an outbreak of  war or a new discovery may
cause the price of a stock to make a sudden jump. 
Since the operators
corresponding to jump processes are non-local, one would like to 
consider operators that are the sum of  
 an elliptic operator and a non-local term.

Such operators  are not yet well understood.
In order to study them  and the
influence of the non-local part, it is quite natural to first look
at the extreme case, that is, where the operator has no differential part, and to begin by  understanding the 
potential theory, existence and uniqueness questions, and 
stochastic differential equations  for 
non-local operators and the associated pure jump processes.

The first such purely non-local operator one would want to study is
the fractional Laplacian $-(-\Delta)^{\al/2}$, where $\Delta$ is the
Laplacian and $\al\in (0,2)$. Such operators have been much studied;
the stochastic processes associated
to these operators are known as symmetric stable processes. See \cite{bogdan1},
\cite{bogdan2},  and \cite{bogdan3}
for a sampling of research on these processes and operators.

The next simplest class of operators $\sL$ is a class introduced
in \cite{jpharn}, known as  stable-like
operators. These are operators $\sL$  defined by
\bee\label{I-def1}
\sL f(x)=\int_{\R^d\setminus\{0\}}
[f(x+h)-f(x)-1_{(|h|\leq 1)}\grad f(x)\cdot h] \frac{A(x,h)}{|h|^{d+\al}}\, dh
\eee
for $f\in C^2(\R^d)$
when $\al\in [1,2)$ and 
\bee\label{I-def2}
\int_{\R^d\setminus\{0\}} [f(x+h)-f(x)] \frac{A(x,h)}{|h|^{d+\al}}\, dh
\eee when $\al\in (0,1)$.
We use $x\cdot y$ for the inner product in $\R^d$.
These stable-like operators
 bear the same relationship to the fractional Laplacian as elliptic
operators in non-divergence form do to the usual Laplacian. 
The name stable-like (which was introduced in \cite{upjump} and
also used in \cite{CK03}) refers to the fact that the jump intensity measure
$A(x,h)/|h|^{d+\al}\, dh$ is comparable to that of the jump intensity measure
of a symmetric stable process.
See \cite{jpharn}, \cite{huili}, \cite{lim},  and \cite{song-vondracek} 
 for some additional results on these operators.
See \cite{nlharn}, \cite{mosernl}, \cite{conhar}, \cite{bkk}, \cite{hadj},
\cite{CK03}, \cite{CK08}, \cite{foondun1}, \cite{foondun2}, 
\cite{kassmann-thesis}, \cite{kolo}, and \cite{huili-thesis}
 for results on operators that are very closely related  to \eqref{I-def1}
and \eqref{I-def2} and
which are also sometimes known as stable-like operators.

Two of the first questions one might ask about stable-like operators
given by \eqref{I-def1} and \eqref{I-def2}  are the H\"older continuity  of harmonic functions and whether a Harnack inequality
holds for non-negative functions  that are harmonic with respect to $\sL$
when the function $A(x,h)$ only satisfies some boundedness and measurability
conditions. These questions were answered in \cite{jpharn}; see also 
\cite{kassmann-thesis} and \cite{song-vondracek}.
A natural  question one might then  ask is whether one can assert additional smoothness
for the solution $u$ to the equation $\sL u=f$ if $A(x,h)$ and $f$  also satisfy some
continuity conditions. The answer to this last
question is the subject of this paper.

Let $\al\in (0,2)$. We impose the following conditions on $A(x,h)$.

\begin{assumption}\label{A1}
Suppose
\begin{enumerate}
\item There exist positive finite constants $c_1,c_2$ such that
$$c_1\leq A(x,h)\leq c_2, \qq x,h\in \R^d.$$
\item There exist $\beta\in (0,1)$ and 
a positive constant $c_3$ such that
$$\sup_{x} \sup_h |A(x+k,h)-A(x,h)|\leq c_3|k|^\beta, \qq k\in \R^d.$$
\item Neither $\beta$ nor $\al+\beta$ is an integer.
\end{enumerate}
\end{assumption}

The assumption that $A(x,h)$ is uniformly bounded above and below is
the analog of strict ellipticity for an elliptic operator in non-divergence
form. The uniform H\"older continuity of $A(x,h)$  in $x$ is the analog of the 
usual assumptions of H\"older continuity in the Schauder theory; see \cite[Chapter 6]{GT}.
Note that no continuity in $h$ is required here.
Finally, the requirement that neither $\beta$ nor $\apb$ be an integer
is quite reasonable; in the theory of elliptic operators, most estimates
break down when the coefficients are not in a H\"older space of non-integer
order.

Our main result is the following. We let $C^\beta$ and $C^{\apb}$
be the usual H\"older spaces. (We recall the definition in \eqref{HSE1}.)

\bet\label{main}
Let $\sL$ be given by \eqref{I-def1} or \eqref{I-def2}  and suppose Assumptions \ref{A1} hold.
If $u\in C^\apb(\R^d)$ satisfies $\sL u=f$, then the following a priori estimate holds:
there exists $c_1$ not depending on $f$ such that
\bee\label{mainest}
\cabnorm{u}\leq c_1\linorm{u}+c_1\cbetanorm{f}.
\eee
\eet

\ni This is the exact analog of the corresponding estimate for
elliptic operators; see \cite[Chapter 6]{GT}.

Lim \cite{lim} has obtained some partial results along the lines of
Theorem \ref{main}. Our result here extends his results by weakening
the hypotheses and strengthening the conclusions. We
show in Section \ref{S:FRR} that our result is sharp in several respects.

Two additional motivations for Theorem \ref{main} are the following.
In \cite{jpharn} harmonic functions for $\sL$ were discussed. There a probabilistic definition of harmonic functions was 
given because in general a harmonic
function, although H\"older continuous, will not be smooth enough to
be in the domain of $\sL$. 
This is not surprising, because for elliptic operators this is
also the case. 
Theorem \ref{main} gives a sufficient condition
for the harmonic function to be in the domain of $\sL$.
Secondly, 
when one considers the process associated
with $\sL$, an essential tool is, as might
be expected, Ito's formula. However the hypotheses of Ito's formula
require the function to be $C^2$. Therefore it would be useful to have
conditions under which a class of functions associated with the process
are at least $C^2$.

Our proof follows roughly along the lines of the Schauder theory for
elliptic equation as presented in \cite[Chapter 6]{GT}.
There are some major differences, however. The estimates for the 
case when $A(x,h)$ is constant in $x$
are much more difficult than the corresponding estimates for the
Laplacian. In addition, because we are dealing with non-local operators, our 
localization procedure is necessarily quite different.

In Section \ref{S:prelim} we define the H\"older spaces and prove a few estimates
that we will need. Section \ref{S:Dpt} investigates the derivatives of
the semigroup corresponding to the operator $\sL$ in the case when 
$A(x,h)$ does not depend on $x$, while Section \ref{S:resolv} is
concerned with the smoothing properties of the corresponding potential
operator. In Section \ref{S:fdd} we obtain estimates on the integrands
in \eqref{I-def1} and \eqref{I-def2}, and we prove Theorem \ref{main}
in Section \ref{S:F}.

We prove a number of results related to Theorem \ref{main} in Section
\ref{S:FRR}. For example we examine what happens when we add to $\sL$
a zero order term or a first order differential term and what happens when
$A(x,h)$ has further smoothness in $x$. We also discuss there a number of
directions for further research, including 
the 
Dirichlet problem for bounded domains, 
boundary estimates for bounded domains, the parabolic case, the symmetric
jump process case, and the case of variable order operators.

The letter $c$ with subscripts denotes a finite positive constant whose
value may vary from place to place.

\section{H\"older spaces}\label{S:prelim}

Let $\beta\in (0,1)$.
We define the seminorm
\bee\label{HE102}
[f]_{C^\beta} =\sup_{x\in \R^d} \sup_{|h|>0}  \frac{|f(x+h)-f(x)|}
{|h|^\beta}
\eee
and the norm

\bee\label{HE103}
\cbetanorm{f}=\linorm{f}+[f]_{C^\beta},
\eee
and say 
$f$ is H\"older continuous of order $\beta$
if $\cbetanorm{f}<\infty$.

We write $D_if$ for ${\del f}/{\del x_i}$, $D_{ij}f$ for $\del^2f/\del x_i \del x_j$, and so on.
Suppose $\beta>1$ is not an integer and let $m$ be the largest integer strictly less than $\beta$.
We define
\bee\label{HSE1}
\cbetanorm{f}=\linorm{f}+\sum_{j_1, \ldots, j_m=1}^d [D_{j_1\cdots j_m}f]_\beta
\eee 
and say $f\in C^\beta$ if $\cbetanorm{f}<\infty$. It
is well known (see the proof of Proposition \ref{HP1} below, for example) that this
norm is equivalent to the norm
\begin{align}
\linorm{f}&+\sum_{j_1=1}^d \linorm{D_{j_1}f}+\sum_{j_1,j_2=1}^d \linorm{D_{j_1j_2}f} +\cdots
+\sum_{j_1, \ldots, j_m=1}^d \linorm{D_{j_1\cdots j_m}f}\nn\\
&+\sum_{j_1, \ldots, j_m=1}^d [D_{j_1\cdots j_m}f]_\beta. \label{HSE2}
\end{align}
(When we say two norms $\norm{\cdot}_1$ and $\norm{\cdot}_2$
are equivalent, we mean that there exist constants $c_1, c_2$ such
that
$$c_1\norm{f}_1\leq \norm{f}_2\leq c_2\norm{f}_1$$
for all $f$.)
       
We also use the fact that the $C^\beta$ norm is equivalent to a
second difference norm:
by \cite[Proposition 8 of Chapter V]{stein}, we have

\bep\label{Pr-P1}
For $\beta\in (0,1)\cup (1,2)$, $f\in C^\beta$ if and only if
$f\in L^\infty$ and
there exists $c_1$ such that
$$|f(x+h)+f(x-h)-2f(x)|\leq c_1|h|^\beta, \qq h,x\in \R^d.$$  
The norm
\bee\label{rhc-E28}
\linorm{f}+\sup_x \sup_{|h|>0}  \frac{|f(x+h)+f(x-h)-2f(x)|}{|h|^\beta}
\eee
is equivalent to the $C^\beta$ norm.
\eep

We will sometimes use the notation
$$\linorm{Df}=\sum_{i=1}^d \linorm{D_if}, \qq 
\linorm{D^2f}=\sum_{i,j=1}^d \linorm{D_{ij}f}.$$
In order to be able to include the case of integer $\beta$ in the
next two results, we introduce the
following notation. 
If $a$ is not an integer,  set $N(f,a)=\canorm{f}$;
if $a=1$,  set $N(f,a)=\linorm{f}+\linorm{Df}$; and if $a=2$, set
$N(f,a)=\linorm{f}+\linorm{Df}+\linorm{D^2f}$.
The following proposition is similar to known results.

\bep\label{HP1} 
If $0<a<b<3$ and $\eps>0$, there exists $c_1$ depending 
only on $a, b$, and $\eps$
such that
\bee\label{HE1}
N(f,a)\leq c_1\linorm{f}+\eps N(f,b).    
\eee
\eep

\proof We first do the case when $0<a<b\leq 1$. Let $h_0=\eps^{1/(b-a)}$.
If $|h|<h_0$, then
$$|f(x+h)-f(x)|\leq N(f,b)|h|^b< N(f,b)|h|^a \eps.$$
If $|h|\geq h_0$, then 
$$|f(x+h)-f(x)|\leq \frac{2}{h_0^a} \linorm{f} |h|^a.$$
Combining,   we have
$$\sup_{|h|>0} \frac{|f(x+h)-f(x)|}{|h|^a}\leq \eps N(f,b)+c_2
\linorm{f}.$$
Taking the supremum over $x$, \eqref{HE1} follows immediately.

Second, we do the case $a=1$ and $b\in (1,2]$. Fix $1\leq i\leq d$ and
let $x_0$ be a point in $\R^d$. The case when $N(f,b)=0$ is trivial,
so we suppose not. Let $R=(\linorm{f}/N(f,b))^{1/b}$. 
By the mean value theorem, there exists $x'$ on the line segment between $x_0$ and $x_0+Re_i$
such that
$$|D_if(x')| =\frac{|f(x_0+Re_i)-f(x_0)|}{R}\leq 2\frac{\linorm{f}}{R}.$$
Then 
$$|D_if(x_0)|\leq |D_if(x')|+|D_if(x')-D_if(x_0)|
\leq \frac{2\linorm{f}}{R}+N(f,b)R^{b-1}.$$
With our choice of $R$,
\bee\label{HE2}
 |D_if(x_0)|\leq c_3\linorm{f}^{1-1/b}N(f,b)^{1/b}.
\eee
Taking the supremum over $x_0\in \R^d$ and then applying the inequality
\bee\label{HE21}
x^\theta y^{1-\theta}\leq x+y, \qq  x, y>0, \quad \theta\in (0,1),
\eee
 we
obtain 
$$\linorm{D_if}\leq \frac{c_4}{\eps}\linorm{f}+ \eps N(f,b).$$

Third, suppose $a=2$ and $b\in (2,3)$.  
Applying \eqref{HE2}
with $f$ replaced by $D_{j}f$ and $b$ replaced by $b-1$ and setting $\gamma
=1/(b-1)$, we have
$$\linorm{D_{ij}f}\leq  c_3\linorm{D_jf}^{1- \gamma}\norm{D_jf}_{C^{b-1}}^\gamma.$$
Using the well known inequality $\linorm{g'}\leq c_{5}\linorm{g}^{1/2}\linorm{g''}^{1/2}$
(this is a special case of  \eqref{HE2}) and summing over $i$ and $j$, we have
$$\linorm{D^2f}\leq c_6\linorm{f}^{(1-\gamma)/2}
\linorm{D^2f}^{(1-\gamma)/2}\cbnorm{f}^{\gamma},$$
and therefore
$$\linorm{D^2f} \leq c_7 \linorm{f}^{(1-\gamma)/(1+\gamma)}
\cbnorm{f}^{2\gamma/(1+\gamma)}.$$
Applying \eqref{HE21} with $\theta=(1-\gamma)/(1+\gamma)$,
we obtain \eqref{HE1}.

For the case $a\in (0,1]$ and $b\in (1,2]$, using the
first and second cases above we have
$$N(f,a)\leq c_8\linorm{f}+c_8\linorm{Df}\leq c_8\linorm{f}+c_9\linorm{f}
+\eps N(f,b),$$
and  the remaining cases are
treated similarly.
\qed

\bel\label{HL1}
If $a\in (0,3)$, there exists $c_1$ such that 
$$N(fg,a)  \leq c_1 N(f,a)N(g,a).$$
\eel

\proof Clearly $\linorm{fg}\leq \linorm{f}\linorm{g}$. If $a\in (0,1)$, we
write
$$f(x+h)g(x+h)-f(x)g(x)=f(x+h)[g(x+h)-g(x)]+ g(x)[f(x+h)-f(x)],$$
and it follows that
$$[fg]_{C^a}\leq \linorm{f}\canorm{g}+\linorm{g}\cbnorm{g}.$$

If $a\in (1,2)$, we use $D_i(fg)=f(D_ig)+(D_if)g$. 
As in the above paragraph, we bound
$$[(D_if)g]_{C^{a-1}}\leq \linorm{D_if}\norm{g}_{C^{a-1}}
+\linorm{g}\norm{D_if}_{C^{a-1}}\leq c_2\canorm{f}\canorm{g},$$
and we bound $[f(D_ig)]_{C^a}$ similarly. Doing this for each $i$ takes
care of the case $a\in (1,2)$. 

Similarly, if $a\in (2,3)$, we use 
\bee\label{HE101}
D_{ij}(fg)=f(D_{ij}g)+g(D_{ij}f)
+(D_if)(D_jg)+(D_jf)(D_ig)
\eee As in the first paragraph,
$$[(D_if)(D_jg)]_{C^{a-2}}\leq c_3\norm{D_if}_{C^{a-2}}\norm{D_jg}_{C^{a-2}}
\leq c_4\canorm{f}\canorm{g},$$
The other terms in \eqref{HE101} are similar.

The remaining cases, when $a=1$ and $a=2$, are easy and are left to the
reader.
\qed

We will need the   following lemma.

\bel\label{rhc-L1}
Let $\beta\in (0,1)$.
Let $\vp$ be a nonnegative $C^\infty$ symmetric function with compact support
such that $\int \vp(x)\, dx=1$, and let $\vp_\eps(x)=\eps^{-d} \vp(x/\eps)$.
Define $f_\eps=f*\vp_\eps$. Then
there exists $c_1$ such that 
for each $i$ and $j$
\begin{align}
\linorm{f-f_\eps}&\leq c_1\cbetanorm{f} \eps^\beta,\label{rhc-E61}\\
\linorm{D_i f_\eps} &\leq c_1\cbetanorm{f} \eps^{\beta-1},\qq\mbox{ and}
\label{rhc-E62}\\
\linorm{D_{ij} f_\eps}& \leq c_1 \cbetanorm{f}\eps^{\beta-2}.\label{rhc-E63}
\end{align}
\eel

\proof The first inequality follows from 
\begin{align*}
|f(x)-f_\eps(x)|&=\Big|\int [f(x)-f(x-y)]\vp_\eps(y)\, dy\Big|\\
&\leq \cbetanorm{f}\int |y|^\beta \vp_\eps(y)\, dy\\
&=c_2\cbetanorm{f} \eps^\beta.
\end{align*}
Since $\int D_i \vp_\eps(y)\, dy=0$,
\begin{align*}|D_i f_\eps(x)|
&=\Big|\int [f(x-y)-f(x)]D_i\vp_\eps(y)\, dy\Big|\\
&\leq \cbetanorm{f}\int |y|^\beta |D_i \vp_\eps(y)|\, dy\\
&= c_3\cbetanorm{f} \eps^{\beta-1}.
\end{align*}

Similarly, since $\int D_{ij}\vp_\eps(y)\, dy=0$, then
\begin{align*}
|D_{ij}f_\eps(x)|&=\Big|\int [f(x-y)-f(x)]D_{ij} \vp_\eps(y)\, dy\Big|\\
&\leq \cbetanorm{f}\int |y|^\beta |D_{ij}\vp_\eps(y)|\, dy\\
&=c_4\cbetanorm{f}\eps^{\beta-2}.
\end{align*}
\qed

\section{Derivatives of semigroups}\label{S:Dpt}

Let $Q_t$ be the semigroup of a symmetric stable process
of order $\al$
and let
$q(t,x)$ be the density, that is, 
the function such that $Q_tf(x)=\int f(y)q(t,x-y)\, dy$.
It is well known that $q$ can be taken to be $C^\infty$ in $x$.

\bep\label{Dpt-P21} For each $k>0$ and each $j_1, \ldots, j_k=1, \ldots, d$,
 we have
$$\int |D_{j_1\cdots j_k} q(1,x)|\, dx<\infty.$$
\eep

This can be proved by generalizing the ideas of \cite[Proposition 2.6]{kolo},
which considers the case of first derivatives. 
See also \cite{huili-thesis}. It can also be proved using Fourier transforms
and complex analytic techniques; see
\cite{sztonyk}, for example. We give a simple proof based on subordination.

\proof Let $W_t$ be a $d$-dimensional Brownian motion and let
$T_t$ be a one-dimensional one-sided stable process of index $\al/2$
independent of $W$.
Then it is well known, by the principle of subordination \cite[Section X.7]{Feller}, that
$X_t=W_{T_t}$ is a symmetric stable process of index $\al$. 
Hence
\bee\label{Dpt-E21}
q(1,x)=\int_0^\infty r(t,x)\,\P(T_1\in dt),
\eee
where $r(t,x)=(2\pi t)^{-d/2} e^{-|x|^2/2t}$ is the density of $W_t$.

The number of jumps of $T_t$ of size larger than $\lam$ is a Poisson process
with parameter $c_1 \lam^{-\al/2}$. So the probability that $T_t$ has no
jumps of size $\lam$ or larger by time $1$ is bounded by $\exp(-c_1\lam^{-\al/2})$. 
Because $T_t$ is non-decreasing, 
this implies
$$\P(T_1\leq \lam)\leq \exp(-c_1 \lam^{-\al/2}).$$
Hence for any $N>0$,
\begin{align}
\int_0^\infty (1+t^{-N}) \,\P(T_1\in dt)&\leq 2+\int_0^1 t^{-N}\,\P(T_1\in dt)\label{Dpt-E22}\\
&\leq 2+\sum_{n=0}^\infty 2^{N(n+1)} \,\P(T_1\in [2^{-n-1}, 2^{-n}])\nn\\
&\leq 2+\sum_{n=0}^\infty 2^{N(n+1)} \,\P(T_1\leq 2^{-n})\nn\\
&\leq 2+\sum_{n=0}^\infty 2^{N(n+1)} e^{-c_1(2^{-n})^{\al/2}}<\infty.\nn
\end{align}
It is easy to see that for each $a>0$ there exist $b$ and $c_2$ 
depending on $a$ such that
$$\sup_x (1+|x|^a) r(t,x)\leq c_2(1+t^{-b}), \qq t>0.$$
This and \eqref{Dpt-E22} allow us to use dominated convergence to differentiate under the
integral sign in \eqref{Dpt-E21}, and we obtain
$$D_{j_1 \cdots j_k} q(1,x)=\int_0^\infty D_{j_1 \cdots j_k} r(t,x)\,\P(T_1\in dt).$$
Then, using \eqref{Dpt-E22} again and Fubini,
\begin{align*}
\int |D_{j_1 \cdots j_k} q(1,x)|\, dx&\leq \int_0^\infty \int |D_{j_1 \cdots j_k} r(t,x)|\, dx\,\P(T_1\in dt)\\
&\leq c_3\int_0^\infty t^{-k /2} \,\P(T_1\in dt)<\infty.
\end{align*}
\qed

If $f\in L^\infty$, it follows easily that $Q_1 f$ is $C^\infty$ for
$t>0$ and for each $j_1, \ldots, j_k$
$$|D_{j_1 \cdots j_k} Q_1 f(x)|\leq c_1\linorm{f}.$$
By scaling we have
\bee\label{Dpt-E23}
|D_{j_1 \cdots j_k} Q_tf(x)|\leq ct^{-k/\al}\linorm{f}.
\eee

Now we consider L\'evy processes whose L\'evy measure is comparable
to that of a symmetric stable process of index $\al$. Suppose
$A_0:\R^d\setminus\{0\}\to [\kappa_1, \kappa_2]$, where $\kappa_1, \kappa_2$
are finite positive constants. Define
\bee\label{Dpt-E10}
\sL_0 f(x)=\int_{\R^d\setminus\{0\}} [f(x+h)-f(x)-1_{(|h|\leq 1)} \grad f(x)\cdot h]
\frac{A_0(h)}{|h|^{d+\al}}\, dh
\eee
for $C^2$ functions $f$ when $\al\geq 1$, and without the
$\grad f(x)$ term when $\al<1$.  Let $P_t$ be the semigroup corresponding to
 the 
 generator $\sL_0$.

\bet\label{Dpt-T1}
If $f\in L^\infty$,
then $P_tf$ is $C^\infty$
for $t>0$ and for each $j_1, \ldots, j_k=1, \ldots, d$,
there exists $c_1$  (depending on $k$) such that
$$|D_{j_1 \cdots j_k}  P_tf(x)|\leq c_1 t^{-k/\al} \linorm{f}.$$
\eet

\proof Let $\sL_1$ be defined by \eqref{Dpt-E10} but with
$A_0(h)$ replaced by $\kappa_1$ and let $\sL_2=\sL_0-\sL_1$. Let
$Q^1_t$ and $Q^2_t$ be the semigroups for the L\'evy processes
with generators $\sL_1, \sL_2$, resp., and let $X^1$, $X^2$ be the
corresponding L\'evy processes. If we take $X^1$ independent of $X^2$, then 
$X^1+X^2$ has the law of the L\'evy process corresponding to
the generator $\sL$. Therefore $P_t=Q^2_tQ_t^1$. We know that
$Q^1_tf$ satisfies the desired estimate by \eqref{Dpt-E23}
and the fact that the process associated with $\sL_1$ is a deterministic
time change of the process considered in Proposition \ref{Dpt-P21}.
By translation invariance, $Q^2_t$ commutes with differentiation.
Therefore  $P_tf=Q_t^2Q_t^1f$ also satisfies the desired estimate,
since
$$\linorm{D_{j_1 \cdots j_k} P_tf}=\linorm{Q_t^2 D_{j_1 \cdots j_k} Q^1_tf}\leq \linorm{D_{j_1 \cdots j_k} Q^1_tf}\leq
c_1 t^{-k/\al}\linorm{f}.$$
\qed

\section{Potentials and H\"older continuity}\label{S:resolv}

Let $P_t$ continue to be the semigroup corresponding to the L\'evy process
in $\R^d$ 
with infinitesimal generator $\sL_0$ given by \eqref{Dpt-E10}
and define the potential
$$Rf(x)=\int_0^\infty P_tf(x)\, dt$$
when the function $t\to P_tf(x)$ is integrable.
We want to prove that $R$ takes functions in $C^\beta$  into
functions in $C^{\al+\beta}$, provided neither $\beta$ nor $\al+\beta$ is
an integer and that $Rf$ is bounded.

\bep\label{rhc-P1} Suppose $\beta\in (0,1)$, $f\in C^\beta$, $Rf\in L^\infty$,
 and $\al+\beta<1$. Then 
$Rf\in C^{\al+\beta}$ and
there exists $c_1$  not depending on $f$
such that  
$\cabnorm{Rf}\leq c_1\cbetanorm{f}+c_1\linorm{Rf}.$
\eep

\proof We first prove that
\bee\label{rhc-E23}
|P_sf(x)-P_sf(y)|\leq \frac{c_2}{s^{(1-\beta)/\al}}|y-x|\, \cbetanorm{f}.
\eee
Define $f_\eps$ as in Lemma \ref{rhc-L1}.

We have, using Theorem \ref{Dpt-T1} and \eqref{rhc-E61},
\begin{align}
|P_s(f-f_\eps)(y)-P_s(f-f_\eps)(x)|&\leq \linorm{\grad P_s(f-f_\eps)} |y-x|\label{rhc-E13}\\
&\leq \frac{c_3}{s^{1/\al}}\linorm{f-f_\eps} |y-x|\nn\\
&\leq \frac{c_3}{s^{1/\al}}\eps^\beta \cbetanorm{f}|y-x|.\nn
\end{align}
Also, using \eqref{rhc-E62},
\begin{align}
|P_sf_\eps(y)-P_sf_\eps(x)|&\leq c_4\linorm{\grad P_s f_\eps} |y-x|\label{rhc-E14}\\
&\leq c_4\linorm{\grad f_\eps} |y-x|\nn\\
&\leq c_5\eps^{\beta-1}\cbetanorm{f} |y-x|. \nn
\end{align}
Setting $\eps=s^{1/\al}$ and combining \eqref{rhc-E13} and \eqref{rhc-E14}
yields \eqref{rhc-E23}.

If $x,y\in \R^d$ and we define $g(z)=f(y-x+z)$, then by the translation
in\-var\-iance of $P_s$ (that is, $P_s$ commutes with translation), $P_sg(x)=P_sf(y)$, and then
$$|P_sf(y)-P_sf(x)|=|P_s(g-f)(x)|\leq \linorm{g-f}\leq \cbetanorm{f}|y-x|^\beta.$$
So putting $t_0=|y-x|^\al$,  we have
\bee\label{PH-E1}
\int_0^{t_0} |P_sf(y)-P_sf(x)|\, ds\leq t_0\cbetanorm{f}|y-x|^\beta=\cbetanorm{f} |y-x|^\apb.
\eee
Using \eqref{rhc-E23} and noting $(1-\beta)/\al>1$,
\begin{align*}
\int_{t_0}^\infty |P_sf(y)-P_sf(x)|\, ds&\leq \int_{t_0}^\infty
\frac{c_6}{s^{(1-\beta)/\al}}\cbetanorm{f}|y-x|\, ds\\
&=c_7t_0^{1-(1-\beta)/\al}\cbetanorm{f}|y-x|\\
&=c_7\cbetanorm{f}|y-x|^{\al+\beta}.
\end{align*}
Combining this with \eqref{PH-E1} and the fact that
\bee\label{rhc-E15}
|Rf(y)-Rf(x)|\leq \int_0^{t_0} |P_sf(y)-P_sf(x)|\, ds
+\int_{t_0}^\infty |P_sf(y)-P_sf(x)|\, ds,
\eee
our result follows.
\qed

Next we consider the case when $0<\beta<1$ and $1<\al+\beta<2$.

\bep\label{rhc-P2} Suppose $\beta\in (0,1)$, $f\in C^\beta$, $\linorm{Rf}<\infty$,
 and $\al+\beta\in (1,2)$. Then 
$Rf\in C^{\al+\beta}$ and
there exists $c_1$ not depending
on $f$ such that 
$\cabnorm{Rf}\leq c_1\cbetanorm{f}+c_1\linorm{Rf}.$
\eep

\proof
Define
$$V_{hs}(f)(x)=P_sf(x+h)+P_sf(x-h)-2P_sf(x).$$
First we show
\bee\label{rhc-E41}
|V_{hs}(f)(x)|\leq c_2|h|^2 \cbetanorm{f} s^{-(2-\beta)/\al}.
\eee
By  Theorem \ref{Dpt-T1}, \eqref{rhc-E61}, and Taylor's theorem,
\begin{align}
|V_{hs}(f-f_\eps)(x)|&\leq c_3|h|^2\linorm{D^2P_s(f-f_\eps)}\label{rhc-E42}\\
&\leq \frac{c_4}{s^{2/\al}}|h|^2 \linorm{f-f_\eps}\nn\\
&\leq \frac{c_5}{s^{2/\al}}|h|^2 \eps^\beta \cbetanorm{f}.\nn
\end{align}
If we set $g_{1\eps}(z)=f_\eps(z+h)$ and $g_{2\eps}(z)=f_\eps(z-h)$, by 
the translation invariance of $P_s$,
$V_{hs}(f_\eps)(x)=P_sg_{1\eps}(x)+P_sg_{2\eps}(x)-2P_sf_\eps(x)$,
and therefore
by \eqref{rhc-E63}
\begin{align}
|V_{hs}(f_\eps)(x)|&=|P_s(g_{1\eps}+g_{2\eps}-2f_\eps)(x)|
\leq \linorm{g_{1\eps}+g_{2\eps}-2f_\eps}\label{rhc-E425}\\
&\leq c_6|h|^2\, \linorm{D^2f_\eps} 
\leq c_7|h|^2\eps^{\beta-2} \cbetanorm{f}.\nn
\end{align}
Letting $\eps=s^{1/\al}$ and combining with \eqref{rhc-E42}, we obtain \eqref{rhc-E41}.

Using \eqref{rhc-E41} and noting $(2-\beta)/\al>1$,
\bee\label{rhc-E44}
\int_{|h|^\al}^\infty |V_{hs}(f)(x)|\leq c_8|h|^2\cbetanorm{f}\int_{|h|^\al}
^\infty s^{-(2-\beta)/\al}\, ds=c_9\cbetanorm{f} |h|^{\al+\beta}.
\eee
Let $g_{10}(x)=f(x+h)$, $g_{20}(x)=f(x-h)$. By translation invariance
and
the H\"older continuity of $f$, 
\begin{align*}
|V_{hs}(f)(x)|&=|P_s(g_{10}+g_{20}-2f)(x)|\leq
\linorm{g_{10}+g_{20}-2f}\\
& \leq 2\cbetanorm{f}|h|^\beta,
\end{align*}
 and thus
\bee\label{rhc-E45}
\int_0^{|h|^\al} |V_{hs}(f)(x)|\, ds\leq 2\cbetanorm{f}|h|^{\al+\beta}.
\eee
Adding \eqref{rhc-E44} and \eqref{rhc-E45} we conclude
$$|Rf(x+h)+Rf(x-h)-2Rf(x)|\leq c\cbetanorm{f}|h|^{\al+\beta}.$$
This with Proposition \ref{Pr-P1} completes the proof.
\qed

Finally we consider the case when $\al+\beta\in (2,3)$.

\bep\label{rhc-P3} Suppose $\beta\in (0,1)$, $f\in C^\beta$, 
$\linorm{Rf}<\infty$,
 and $\al+\beta\in (2,3)$. Then $Rf\in C^{\al+\beta}$ and
there exists $c_1$ not depending on $f$ such that
$\cabnorm{Rf}\leq c_1\cbetanorm{f}+c_1\linorm{Rf}.$
\eep

\proof Necessarily $\al>1$. In view of Proposition \ref{HP1} it suffices to show
\bee\label{rhc-E71}
\norm{D_iRf}_{C^{\al+\beta-1}}\leq c_2\cbetanorm{f}, \qq i=1, \dots, d.
\eee
Fix $i$ and let $Q_t=D_iP_t$. 
From Theorem \ref{Dpt-T1} we  have 
$$\linorm{D_{j_1 j_2} Q_tf}
\leq c_3t^{-3/\al}\linorm{f}, \qq j_1, j_2=1, \ldots d.$$ 

Define $W_{hs}(f)(x)=Q_sf(x+h)+Q_sf(x-h)-2Q_sf(x)$. Note that
$Q_s$ is translation invariant. Analogously to
\eqref{rhc-E42} and \eqref{rhc-E425},
\begin{align*}
|W_{hs}(f-f_\eps)(x)|&\leq c_4|h|^2 \linorm{D^2Q_s(f-f_\eps)}\\
&\leq \frac{c_5|h|^2}{s^{3/\al}}\linorm{f-f_\eps}\\
&\leq \frac{c_6|h|^2}{s^{3/\al}}\eps^\beta \cbetanorm{f}
\end{align*}
and
\begin{align*}
|W_{hs}(f_\eps)(x)|
&\leq c_7|h|^2 \linorm{D^2Q_sf_\eps}=c_7|h|^2\linorm{Q_sD^2f_\eps}\\
&\leq \frac{c_8|h|^2}{s^{1/\al}}\linorm{D^2f_\eps}
\leq \frac{c_9}{s^{1/\al}}|h|^2 \eps^{\beta-2}\cbetanorm{f}.
\end{align*}
Taking $\eps=s^{1/\al}$ we obtain
$$|W_{hs}(f)(x)|\leq c_{10}|h|^2 s^{(\beta-3)/\al}\cbetanorm{f}.$$
Integrating this bound   over $[\,|h|^\al,\infty)$ yields
$c_{11}|h|^{\al+\beta-1}\cbetanorm{f}$. 

On the other hand, if $g_{10}$ and $g_{20}$ are defined as in the
proof of Proposition \ref{rhc-P2},
\begin{align*}
|W_{hs}(f)(x)|&\leq \linorm{Q_s(g_{10}+g_{20}-2f)}
\leq c_{12} s^{-1/\al}\linorm{g_{10}+g_{20}-2f}\\
&\leq c_{13} s^{-1/\al}|h|^\beta \cbetanorm{f},
\end{align*}
and integrating this bound over $s$
from 0 to $|h|^\al$ yields $c_{14}|h|^{\al+\beta-1}\cbetanorm{f}$;
we use the fact that $1/\al<1$ here.
Therefore
$$|W_{hs}(D_iRf)(x)|\leq c_{14}|h|^{\al+\beta-1}\cbetanorm{f},$$
which with Proposition \ref{Pr-P1} yields \eqref{rhc-E71}.
\qed

We reformulate and summarize the preceding propositions
in the following theorem.
Let $\sL_0$ be defined as in \eqref{Dpt-E10}.

\bet\label{FTest}
Suppose $\beta\in (0,1)$ and $\apb\in (0,1)\cup (1,2)\cup (2,3)$.
There exists $c_1$ such that if $u$ is in the domain of $\sL_0$ and
$\sL_0 u=f$ with $\linorm{u}<\infty$, then
\bee\label{FEest}
\cabnorm{u}\leq c_1\canorm{f}+c_1\linorm{u}.
\eee
\eet

\proof If $\sL_0u=f$ and $\linorm{u}<\infty$, then we have $u=-Rf$, and so the result follows by Propositions 
\ref{rhc-P1}, \ref{rhc-P2}, and \ref{rhc-P3}.
\qed

\section{First and second differences}\label{S:fdd}

For $f$ bounded define
\bee\label{defEh}
E_hf(x)=f(x+h)-f(x).
\eee
For $f\in C^1$ define
\bee\label{defFh}
F_hf(x)=f(x+h)-f(x)-\grad f(x)\cdot h.
\eee
Observe that if $g:\R\to \R$ is in $C^\gamma$ with $\gamma\in (1,2)$, then
\begin{align}
|g(t)-g(0)-g'(0)t|&=\Big|\int_0^t[g'(s)-g'(0)]\, ds\Big|\label{CE1}\\
&\leq \cgammanorm{g}\int_0^t s^{\gamma-1}\, ds\leq c_1\cgammanorm{g} t^\gamma,\nn
\end{align}
while if $\gamma\in (2,3)$, then
\begin{align}
|g(t)-g(0)-&g'(0)t-\tfrac12 g''(0)t^2|=\Big|\int_0^t [g'(s)-g'(0)]\, ds-\tfrac12 g''(0)t^2\Big|\label{CE2}\\
&=\Big|\int_0^t\int_0^s [g''(r)-g''(0)]\, dr\, ds\Big|\nn\\
&\leq \cgammanorm{g}\int_0^t \int_0^s r^{\gamma-2}\, dr\, ds=c_2\cgammanorm{g}t^\gamma.\nn
\end{align}
Let $Hf$ be the Hessian of $f$, so that 
$$h\cdot Hf(x)k=\sum_{i,j=1}^d h_i D_{ij}f(x)k_j$$
if $h=(h_1, \ldots, h_d)$ and $k=(k_1, \ldots, k_d)$.

\bet\label{fd-T1} Suppose $f\in C^\gamma$ for $\gamma\in (0,1)\cup (1,2)\cup 
(2,3)$. There exists $c_1$ not depending on $f$  such that the following estimates hold.

(a) For all $\gamma$,
\bee\label{oer-E51}
|E_hf(x)|\leq c_1(|h|^{\gamma\land 1}\land 1)
\cgammanorm{f}
\eee
and if $\gamma>1$, 
\bee\label{oer-E52}
|F_hf(x)|\leq c_1(|h|^{\gamma\land 2}\land 1)\cgammanorm{f}.
\eee

(b) For all  $\gamma$, 
\bee\label{oer-E53}
|E_hf(x+k)-E_hf(x)|\leq c_1(|h|^{\gamma\land 1}\land |k|^{\gamma\land 1})\norm{f}_{C^\gamma}.
\eee

(c) If $\gamma\in (1,2)$, then
\bee\label{oer-E54}
|E_hf(x+k)-E_hf(x)|\leq c_1((|h|^{\gamma-1}|k|)\land (|h|\, |k|^{\gamma-1})\norm{f}_{C^\gamma}.
\eee

(d) If $\gamma\in (1,2)$, then
\bee\label{oer-E55}
|F_hf(x+k)-F_hf(x)|\leq c_1((|h|^{\gamma})\land (|h|\, |k|^{\gamma-1}))\norm{f}_{C^\gamma}.
\eee

(e) If $\gamma\in (2,3)$, then 
\bee\label{oer-E57}
|F_hf(x+k)-F_hf(x)|\leq c_1((|k|^{\gamma-2}|h|^2)\land (|h|^{\gamma-1}|k|))\norm{f}_{C^\gamma}.
\eee
\eet

\proof
(a) The estimate for $E_hf$ follows by the definition of $C^\gamma$. 
The one for $F_hf$ follows from \eqref{CE1} or \eqref{CE2} applied
to $g(s)=f(x+sh/|h|)$ with $t=|h|$.

(b) Write
\bee\label{oer-E21A}
E_hf(x+k)-E_hf(x)=[f(x+h+k)-f(x+k)]-[f(x+h)-f(x)],
\eee
and note that because $f\in C^{\gamma}$, this is bounded
by $2|h|^{\gamma\land 1}\norm{f}_{C^{\gamma}}$. We can also write
$E_hf(x+k)-E_hf(x)$ as
\bee\label{oe-E21}[f(x+h+k)-f(x+h)]-[f(x+k)-f(x)],\eee
so we also get the bound $2|k|^{\gamma\land 1}\norm{f}_{C^{\gamma}}$.

(c)  Using \eqref{CE1}
$$f(x+h+k)-f(x+k)=\grad f(x+k)\cdot h+R_1$$
and
$$f(x+h)-f(x)=\grad f(x)\cdot h+R_2,$$
where $R_1$ and $R_2$ are both bounded by $c_2\cgammaf |h|^\gamma$.
By \eqref{oer-E21A}
$$E_hf(x+k)-E_hf(x)=[\grad f(x+k)-\grad f(x)]\cdot h+R_1-R_2,$$
and the right hand side is bounded by
\bee\label{CE11}
c_3\cgammaf (|k|^{\gamma-1}|h|+|h|^\gamma).
\eee
Starting with \eqref{oe-E21} instead of \eqref{oer-E21A} we also
get the bound
\bee\label{CE12}
c_3\cgammaf (|h|^{\gamma-1}|k|+|k|^\gamma).
\eee
Using \eqref{CE11} when $|h|\leq |k|$ and \eqref{CE12} when
$|h|>|k|$ proves \eqref{oer-E54}.
 
(d)
By \eqref{CE2}
$$|F_hf(x)|\leq c_3\cgammaf |h|^\gamma,$$
and the same bound holds for $F_hf(x+k)$,
so 
\bee\label{CE31}
|F_hf(x+k)-F_hf(x)|\leq c_3\cgammaf|h|^\gamma.
\eee
On the other hand
$$f(x+k+h)-f(x+h)=\grad f(x+h)\cdot k+R_3$$
and
$$f(x+k)-f(x)=\grad f(x)\cdot k+R_4,$$
where $R_3$ and $R_4$ are both bounded by $c_4\cgammaf |k|^\gamma$.
Also
$$|\grad f(x+k)\cdot h-\grad f(x)\cdot h|\leq c_5\cgammaf |h|\, |k|^{\gamma-1}$$
and 
$$|\grad f(x+h)\cdot k-\grad f(x)\cdot k|\leq c_5\cgammaf |k|\, |h|^{\gamma-1}.$$
Combining and using the fact that $\gamma<2$,
$$|F_hf(x+k)-F_hf(x)|\leq c_6\cgammaf(|k|^\gamma+|h|\, |k|^{\gamma-1}+|k|\, |h|^{\gamma-1}),$$
which together with \eqref{CE31} proves \eqref{oer-E55}.

(e)
Applying \eqref{CE2}
\bee\label{CE25}
|F_hf(x)-\tfrac12 h\cdot Hf(x)h|\leq c_7 \cgammaf |h|^\gamma
\eee
and we obtain the same bound for $|F_hf(x+k)-\tfrac12 h\cdot Hf(x+k)h|$.
Since
$$|h\cdot (Hf(x+k)-Hf(x))h|\leq c_8\cgammaf |h|^2|k|^{\gamma-2},$$
then
\bee\label{CE3}
|F_hf(x+k)-F_hf(x)|\leq c_9\cgammaf (|h|^2|k|^{\gamma-2}+|h|^\gamma).
\eee
On the other hand, using \eqref{CE1} and \eqref{CE2},
$$f(x+k+h)-f(x+k)=\grad f(x+h)\cdot k+\tfrac12k\cdot Hf(x+h)k+R_5,$$
$$f(x+k)-f(x)=\grad f(x)\cdot k+\tfrac12k\cdot Hf(x)k+R_6,$$
$$\grad f(x+k)\cdot h-\grad f(x)\cdot h=k\cdot Hf(x)h+R_7,$$
and
$$\grad f(x+h)\cdot k-\grad f(x)\cdot k=h\cdot Hf(x)k+R_8,$$
where $R_5$ and $R_6$ are both bounded by $c_{10}\cgammaf |k|^\gamma$,
$R_7$ is bounded by $c_{10}\cgammaf |k|^{\gamma-1}|h|$, and
$R_8$ is bounded $c_{10}\cgammaf |h|^{\gamma-1}|k|$.
Therefore
$$|F_hf(x+k)-F_hf(x) -\tfrac12 k\cdot (Hf(x+h)-Hf(x))k|\leq |R_5|
+|R_6|+|R_7|+|R_8|,$$
which implies
\begin{align}   
|F_hf(x+k)-F_hf(x)|&\leq c_{11}\cgammaf
(|k|^\gamma+|k|^{\gamma-1}|h|\label{CE7}\\
&\qq +|h|^{\gamma-1}|k|+|k|^2|h|^{\gamma-2}).\nn
\end{align}
Using \eqref{CE3} if $|h|\leq |k|$ and \eqref{CE7} if $|h|>|k|$ proves \eqref{oer-E57}.
\qed

We have the following corollary.

\bec\label{fd-C1} Suppose $f\in C^{\al+\beta}$ for some $\beta\in (0,1)$
and $\alpha+\beta\in (0,1)\cup (1,2)\cup (2,3)$. 
There exists $c_1$ not depending on $f$ such that

(a) If $\al<1$, then
\bee\label{fd-EW1}
\int  |E_hf(x+k)-E_hf(x)|\, \frac{dh}{|h|^{d+\al}}
\leq c_1|k|^\beta \cabnorm{f}.
\eee

(b) If $\al\in [1,2)$, then 
\begin{align}
\int_{|h|\leq 1} |F_hf(x+k)-F_hf(x)|\, \frac{dh}{|h|^{d+\al}}
&+\int_{|h|> 1} |E_hf(x+k)-E_hf(x)|\, \frac{dh}{|h|^{d+\al}}\nn\\
&\leq c_1|k|^\beta \cabnorm{f}.\label{fd-EW2}
\end{align}
\eec

\proof 
If $|k|>1$, the left hand side of \eqref{fd-EW1} is less than or
equal to 
$$\int |E_hf(x+k)|\frac{dh}{|h|^{d+\al}}+\int |E_hf(x)|\frac{dh}{|h|^{d+\al}},$$
which is bounded using Theorem \ref{fd-T1}(a).  We treat \eqref{fd-EW2}
similarly.

If $|k|\leq 1$, we use the bounds in Theorem \ref{fd-T1}(b)--(e), breaking the
integrals into three: where $|h|<|k|$, where $|k|\leq |h|\leq 1$, and
where $|h|>1$. The rest is elementary calculus.
\qed

\begin{remark}\label{oe-R111}{\rm
By Theorem \ref{fd-T1}(a), the integrals defining $\sL u$ are thus absolutely convergent
if $u\in C^\apb$ for some $\beta>0$.
In particular, the domain of $\sL$ contains $C^\apb$ for each $\beta>0$.
}
\end{remark}

The following is immediate from Corollary \ref{fd-C1}.
 
\bec\label{fd-C45}
Suppose $u\in C^\apb$ for some $\beta\in (0,1)$ and 
$\apb\in (0,1)\cup (1,2) \cup (2,3)$. Let $\sL_0$ be defined
by \eqref{Dpt-E10}. Then $\sL_0 u\in C^\beta$ and there exists $c_1$
such that
$$\cbetanorm{\sL_0 u}\leq c_1\cabnorm{u}.$$
\eec

\section{Proof of Theorem \ref{main}}\label{S:F}

Let $B(x,r)$ denote the ball of radius $r$ centered at $x$.
 Let $\ol \vp$ be a cut-off function that is  1 on
$B(0,1)$, 0 on $B(0,2)^c$, takes values in $[0,1]$, and is $C^\infty$.
Let $\vp_{r,x_0}(x)=r^{-d}\ol \vp((x-x_0)/r)$. 
When $r$ and $x_0$ are clear, we will write just $\vp$
for $\vp_{r,x_0}$.  

\bep\label{FP1} Suppose $\cabnorm{u}<\infty$. Suppose
for each $\delta>0$ there exists $r$ and $c_1$ (depending on $\delta$)
such that
\bee\label{FE1}
\cabnorm{u\vp_{r,x_0}}\leq c_1\cbetanorm{f}+c_1\linorm{u}
+\delta \cabnorm{u}
\eee
Then there exists $c_2$ depending on $\delta$  such that
\bee\label{FE2}
\cabnorm{u}\leq c_2\cbetanorm{f}+c_2\linorm{u}.
\eee
\eep

\proof
First we do the case where $\al+\beta\in (0,1)\cup (1,2)$.
Recall from Proposition \ref{HP1} that there exist $c_3$ and $c_4$ such
that 
\begin{align}
c_3\cabnorm{g}\leq \linorm{g}&+\sup_x \sup_{|h|>0}
\frac{g(x+h)+g(x-h)-2g(x)}{|h|^{\apb}}\label{FE201}\\
&\leq c_4\cabnorm{g}\nn
\end{align}
for all $g\in C^{\al+\beta}$.
Choose $\delta=c_3/2c_4$ and then choose $r$ and $c_1$ using \eqref{FE1}.
 If $x_0\in \R^d$, let
$v=u\vp_{r,x_0}$, and
note that  $u=v$ in the ball $B(x_0,r)$. 
If $|h|< r$,
\begin{align}
|u(x_0+h)+u(x_0-h)-2u(x_0)|&=|v(x_0+h)+v(x_0-h)-2v(x)|\label{FE11}\\
&\leq c_4\cabnorm{v} |h|^{\al+\beta}.\nn
\end{align}
On the other hand, if $|h|\geq r$,
\bee\label{FE12}
|u(x_0+h)+u(x_0-h)-2u(x_0)|\leq \frac{4}{r^{\al+\beta}}
\linorm{u} |h|^{\al+\beta}=c_5\linorm{u} |h|^{\al+\beta}.
\eee
Combining \eqref{FE11} and \eqref{FE12} and using \eqref{FE1},
\begin{align*}
|u(x_0+h)+&u(x_0-h)-2u(x_0)|\leq (c_4\cabnorm{v}+c_5\linorm{u})
|h|^{\al+\beta}\\
&\leq (c_1c_4\cbetanorm{f}+(c_1c_4+c_5)\linorm{u}+c_4\delta\cabnorm{u})|h|^{\al+\beta}.
\end{align*}
This and \eqref{FE201} yield
$$\cabnorm{u}\leq c_6\cbetanorm{f}+c_6\linorm{u}+\tfrac12\cabnorm{u}.$$
Subtracting $\frac12 \cabnorm{u}$ from both sides and multiplying
by 2 gives \eqref{FE1}.

Now we consider the case when $\al+\beta\in (2,3)$. 
Since $u\in C^{\al+\beta}$ if $u\in L^\infty$ and each $D_iu\in C^{\al+\beta-1}$,
by \eqref{HSE1}, \eqref{HSE2}, and Propositions \ref{Pr-P1} and \ref{HP1} there exists $c_7$ such that
$$\cabnorm{u}\leq c_7 \Big(\linorm{u}+
\sup_i \sup_x \sup_{|h|>0} 
\frac{|D_iu(x+h)+D_iu(x-h)-2D_iu(x)|}{|h|^{\al+\beta-1}}\Big).$$
Let $\delta=1/2c_7(1+c_4)$, choose $r$ using \eqref{FE1}, and let 
$v=u\vp_{r,x_0}$.
If $|h|<r$, then
for any $i$,
\begin{align*}
|D_iu(x_0+h)+&D_iu(x_0-h)-2D_iu(x_0)|\\
&=|D_iv(x_0+h)+D_iv(x_0-h)-2D_iv(x_0)|\\
&\leq c_{4}\cabnorm{v}  |h|^{\al+\beta-1}\\
&\leq (c_1c_{4}\cbetanorm{f}+c_{1}c_4\linorm{u}+\delta c_4\cabnorm{u})\, |h|^{\al+\beta-1}.
\end{align*}
On the other hand, if 
$|h|\geq r$, then
\bee\label{perE1}
|D_iu(x_0+h)+D_iu(x_0-h)-2D_iu(x_0)|\leq \frac{4}{r^{\al+\beta-1}}
\linorm{D_iu} |h|^{\al+\beta-1}.
\eee
Choose $\eps=r^{\apb-1}\delta/4$ and then use Proposition \ref{HP1} to see there
exists $c_{8}$ such that
$$\linorm{D_iu}\leq c_{8}\linorm{u}+\eps \cabnorm{u}.$$
Substituting this in \eqref{perE1},
$$|D_iu(x_0+h)+D_iu(x_0-h)-2D_iu(x_0)|\leq (c_{9}\linorm{u}+ \delta \cabnorm{u})|h|^{\apb-1}.$$
Therefore
\begin{align*}
|D_iu(x_0+h)+&D_iu(x_0-h)-2D_iu(x_0)|\\
&\leq (c_{10}\cbetanorm{f}+
c_{10} \linorm{u}+ (1+c_4)\delta\cabnorm{u}) |h|^{\al+\beta-1},
\end{align*}
and hence
$$\cabnorm{u}\leq c_{11}\cbetanorm{f}+c_{11}\linorm{u}+\tfrac12\cabnorm{u}.$$
Subtracting
$\frac12 \cabnorm{u}$ from both sides, and multiplying by
2 yields our result.
\qed

\longproof{of Theorem \ref{main}}
\ni{\it Step 1.} In this step we define  a certain function $F$. Let us suppose for now that $\al<1$, leaving the case $\al\geq 1$ until later.
Fix $\delta>0$ and let $\eps>0$ be chosen later. Let $x_0\in \R^d$
be fixed and choose $r$ such that
$$\sup_{|h|>0}|A(x,h)-A(x_0,h)|<\eps$$
if $|x-x_0|\leq 4r$.
Let $b(x,h)=A(x,h)-A(x_0,h)$, 
$$\sL_0u(x)=\int [u(x+h)-u(x)] \frac{A(x_0,h)}{|h|^{d+\alpha}}\, dh,$$  
and $\sB=\sL-\sL_0$. 
Let $\vp=\vp_{r,x_0}$ be as in the paragraph preceding Proposition \ref{FP1}
and let $v=u\vp$.

We have
\begin{align*}
v(x+h)-v(x)&=u(x)[\vp(x+h)-\vp(x)]+\vp(x)[u(x+h)-u(x)]\\
&\qq +[u(x+h)-u(x)]\, [\vp(x+h)-\vp(x)],
\end{align*}
and therefore
\begin{align*}\sL v(x)&=u(x)\sL\vp(x)+\vp(x)\sL u(x)+H(x)\\
&=u(x)\sL\vp(x)+\vp(x)f(x)+H(x),
\end{align*}
where
$$H(x)=\int [u(x+h)-u(x)]\, [\vp(x+h)-\vp(x)] \frac{A(x,h)}{|h|^{d+\al}}\, dh.$$
On the other hand, 
$$\sL v(x)=\sL_0 v(x)+\sB v(x),$$
and so  we have
\begin{align}   
\sL_0 v(x)&=u(x) \sL \vp(x)+\vp(x) f(x)+H(x)-\sB v(x)\label{FE4}\\
&=J_1(x)+J_2(x)+J_3(x)+J_4(x).\nn
\end{align}
Set 
\bee\label{FQE3}
F(x)=\sum_{i=1}^4 J_i(x).
\eee

By Theorem \ref{FTest}
we have
$$\cabnorm{v}\leq c_1(\cbetanorm{F} + \linorm{v})\leq c_1(\cbetanorm{F}+\linorm{u}).$$
So if, given $\eps$, we can  show
\bee\label{FE3}
\cbetanorm{F}\leq c_2(\cbetanorm{f}+\linorm{u}+\eps
\cabnorm{u}),
\eee
we take $\eps=\delta/c_2$,  we then have \eqref{FE1},  
we apply Proposition \ref{FP1}, and we are done.

\ni{\it Step 2.}
We first look at the $L^\infty$ norm of $F$. Since
$$\int [\vp(x+h)-\vp(x)] \frac{1}{|h|^{d+\al}}\, dh
\leq \int |E_h\vp(x)| \frac{1}{|h|^{d+\al}}\, dh\leq c_3<\infty, $$
where $E_h$ is defined in \eqref{defEh},
then
$$|u(x)\sL \vp(x)|\leq c_3 \linorm{u}.$$
Similarly
\begin{align*}
|H(x)|&=\Big|\int [u(x+h)-u(x)]\, [\vp(x+h)-\vp(x)]
\frac{A(x_0,r)}{|h|^{d+\al}}\, dh\Big|\\
&\leq c_4\linorm{u} \int |E_h\vp(x)| \frac{1}{|h|^{d+\al}}\, dh\\
&\leq c_5 \linorm{u}.
\end{align*}
We also  have
$$|\vp(x)f(x)|\leq \linorm{f}\leq \cbetanorm{f}.$$

It remains to bound $\sB v(x)$.
If $x\notin B(x_0, 3r)$, then since 
$v(x)=0$ and $v(x+h)=0$ unless $|h|>r$, we see
$$|\sB v(x)|=\Big| \int_{|h|>r} v(x+h) \frac{b(x,h)}{|h|^{d+\al}}\, dh\Big|
\leq c_6\linorm{u} \int_{|h|>r} |h|^{-d-\al}\, dh
=c_7\linorm{u}.$$
We have
$$\cabnorm{v}\leq c_8\cabnorm{\vp}\cabnorm{u}
\leq c_9\cabnorm{u},$$
since $\vp$ is smooth. 
By  Theorem \ref{fd-T1}(a), 
$$|E_hv(x)|\leq c_{10}(|h|^{(\apb)\land 1}\land 1) \cabnorm{v},$$
and so
\begin{align*}
|\sB v(x)| &=\Big|\int E_hv(x) \frac{b(x,h)}{|h|^{d+\al}}\, dh\Big|\\
&\leq c_{10}\eps \int (|h|^{(\apb)\land 1}\land 1) \frac{1}{|h|^{d+\al}}\, dh
\, \cabnorm{v}\\
&\leq c_{11} \eps \cabnorm{v}\leq c_{12}\eps \cabnorm{u}.
\end{align*}
We used the fact that we chose $r$ small so that $|b(x,h)|\leq \eps$.
To summarize, in this step we have shown
\bee\label{FQ1}
\linorm{F}\leq c_{13}(\cbetanorm{f}+\linorm{u}+\eps\cabnorm{u}).
\eee

\ni{\it Step 3.} We next estimate $[F]_{C^{\beta}}$.
Since we have
$$|F(x+k)-F(x)|\leq 2\linorm{F}\leq (2^\beta/r^\beta)\linorm{F}|k|^\beta$$
when $|k|\geq r/2$
and we have an upper bound of the correct form for $\linorm{F}$ in
\eqref{FQ1}, to bound $[F]_{C^\beta}$ it suffices to look at 
$F(x+k)-F(x)$ when
$|k|\leq r/2$. We look at the differences for $J_i$ for $i=1,\ldots, 4$.

We look at $J_4$ first, since this is the most difficult one. 
First suppose $x\notin B(x_0,3r)$. Then $v(x+h+k)$, $v(x+h)$, $v(x+k)$,
and $v(x)$ are all zero if $|h|\leq r/2$. So
\begin{align*}
|\sB v(x+k)-&\sB v(x)|\\
&=\Big| \int_{|h|>r/2} [v(x+h+k)b(x+k,h) -v(x+h)b(x,h)]\frac{dh}{|h|^{d+\al}}\Big|\\
&\leq \int_{|h|> r/2} |v(x+h+k)-v(x+h)|\, |b(x+k,h)|\frac{dh}{|h|^{d+\al}}\\
&\qq +\int_{|h|>r/2} |v(x+h)|\, |b(x+k,h)-b(x,h)|\frac{dh}{|h|^{d+\al}}\\
&\leq c_{14}\cbetanorm{v}|k|^\beta \int_{|h|>r/2} \frac{dh}{|h|^{d+\al}}
+c_{11}\linorm{v}|k|^\beta \int_{|h|>r/2} \frac{dh}{|h|^{d+\al}} .
\end{align*} 
Since $\linorm{v}\leq \linorm{u}$ and
\begin{align*}
\cbetanorm{v}&\leq c_{15}\cbetanorm{u}\cbetanorm{\vp}\leq c_{16}\cbetanorm{u}\\
&\leq c_{17}\linorm{u}+\eps\cabnorm{u}
\end{align*} 
by Proposition \ref{HP1},
we have our required estimate when $x\notin B(x_0,3r)$. 

Now suppose $x\in B(x_0,3r)$. Since $|k|\leq r/2$, then $x+k\in B(x_0,4r)$, 
and so $|b(x,h)|\leq \eps$ and $|b(x+k,h)|\leq \eps$ for all $h$.
We write
\begin{align*}
|\sB v(x+k)&-\sB v(x)|\\
&\leq \int |E_hv(x+k)-E_hv(x)|\frac{|b(x+k,h)|}{|h|^{d+\al}}\, dh\\
&\qq +\int_{|h|\leq \zeta} |E_hv(x)|\frac{|b(x+k,h)-b(x,h)|}{|h|^{d+\al}}\, dh\\
&\qq +\int_{|h|> \zeta} |E_hv(x)|\frac{|b(x+k,h)-b(x,h)|}{|h|^{d+\al}}\, dh\\
&=I_1+I_2+I_3,
\end{align*}
where $\zeta$ will be chosen in a moment.
By Theorem \ref{fd-T1},
$$I_1\leq \eps\int(|h|^{(\apb)\land 1}\land |k|^{(\apb)\land 1})\frac{dh}{|h|^{d+\al}}
\leq c_{18}\eps \cabnorm{v}|k|^\beta.$$
Suppose for the moment that $\al+\beta<1$.
For $I_2$ we have 
$$I_2\leq c_{19}\cabnorm{v}\int_{|h|\leq \zeta} (|h|^\apb\land 1)\frac{|k|^\beta }
{|h|^{d+\al}}\, dh\leq \eps\cabnorm{v},$$
provided we take $\zeta$ small; note that the choice of $\zeta$ can be made to
depend
only on $d$, $\al$, $\beta$, and $\eps$. For $I_3$ we
have
$$I_3\leq c_{20}\cbetanorm{v}\int_{|h|>\zeta} (|h|^\beta\land 1)\frac{|k|^\beta}
{|h|^{d+\al}}\, dh\leq c_{21} \cbetanorm{v} |k|^\beta.$$
We now use
$$\cabnorm{v}\leq c_{22}\cabnorm{u}\cabnorm{\vp}$$
and
\begin{align*}
\cbetanorm{v}&\leq c_{23}\cbetanorm{u}\cbetanorm{\vp}\\
&\leq \eps\cabnorm{u}+ c_{24}\linorm{u}.
\end{align*} 
Summing the estimates for $I_1, I_2$, and $I_3$, 
we have the desired bound for $J_4$ when $\al+\beta<1$.  The case $\al+\beta\in (1,2)$ is very 
similar; the details are left to the
reader.

Next we look at $J_1$. Similarly to the estimates for $J_4$,
we see that $\cbetanorm{\sL \vp}\leq c_{25}$.
We then have
$$\cbetanorm{J_2}\leq c_{26}\cbetanorm{u}\cbetanorm{\sL \vp},$$
and then Proposition \ref{HP1} gives our estimate. 

The estimate for $J_2$ is quite easy.  By Lemma \ref{HL1}
$$\cbetanorm{\vp f}\leq c_{26}\cbetanorm{\vp}\cbetanorm{f}
\leq c_{27}\cbetanorm{f}.$$

It remains to handle $J_3$. We have
\begin{align*}
H(x+k)-H(x)&=\int [E_hu(x+k)-E_hu(x)] E_h\vp(x+k)\frac{A(x+k,h)}{|h|^{d+\al}}\, dh\\
&\qq +\int E_h u(x)[E_h\vp(x+k)-E_h\vp(x)] \frac{A(x+k,h)}{|h|^{d+\al}}\,dh\\
&\qq + \int E_hu(x)E_h\vp(x) \frac{A(x+k,h)-A(x,h)}{|h|^{d+\al}}\, dh\\
&=I_4+I_5+I_6.
\end{align*}
By Theorem \ref{fd-T1}
\begin{align*}
|I_4|&\leq c_{28}|k|^\beta \cbetanorm{u}\int (|h|^\beta\land 1) \frac{dh}{|h|^{d+\al}}\, dh\\
&\leq c_{29}|k|^\beta\cbetanorm{u}.
\end{align*}
Also by Theorem \ref{fd-T1}
\begin{align*}
|I_5|&\leq c_{30}\linorm{u}\int (|h|^\beta\land |k|^\beta \land 1) \frac{dh}{|h|^{d+\al}}\, dh\\
&\leq c_{31}\linorm{u}|k|^\beta;
\end{align*}
to get the second inequality we split the integral into $|h|\leq |k|$, $|k|<|h|
\leq 1$, and $|h|>1$.
Using Theorem \ref{fd-T1} a third time
$$|I_6|\leq c_{32}\linorm{u}\int(|h|^\beta\land 1)\frac{|k|^\beta}{|h|^{d+\al}}\, dh\leq c_{33}\linorm{u}|k|^\beta.$$
Combining yields
$$[H]_{C^\beta}\leq c_{34}\cbetanorm{u},$$
and we now apply Proposition \ref{HP1}.

\ni{\it Step 4.} Finally we consider the case $\al\geq 1$. This is
very similar to the $\al<1$ case, but where we replace the use
of $E_hf$ by $F_hf$. We leave the details to the reader.
\qed

\section{Further results and remarks}\label{S:FRR}

\subsection{An extension}\label{SS1}

We remark that the proof of Theorem \ref{main} really only
required  that there exist $c_1$ and $h_0$ such that
$$\sup_x\sup_{|h|\leq h_0} |A(x+k,h)-A(x,h)|\leq c_1|k|^\beta.$$
The observation  needed is that  one can bound 
$$\Big\Vert \int_{|h|>h_0} [u(x+h)-u(x)] \frac{A(x,h)}{|h|^{d+\al}}\, dh\Big\Vert_{C^\beta}
\leq c_2\cbetanorm{u}\leq      
c_3\linorm{u}+\eps\cabnorm{u}.$$

\subsection{Zero order terms}\label{SS2}

We can add a zero order term to $\sL$ and have the result remain
valid. 

\bet\label{FRRT1}
Let $P$ be a function such that $\cbetanorm{P}<\infty$. Let
$$ \sL' u(x)=\sL u(x)+ P(x) u(x),$$
where $\sL$ is defined by \eqref{I-def1} or \eqref{I-def2} and satisfies the
assumptions of Theorem \ref{main}.
Then there exists $c_1$ (which depends on $\cbetanorm{P}$) such that
if $\sL' u(x)=f(x)$ and $\cabnorm{u}<\infty$, then
$$\cabnorm{u}\leq c_1(\linorm{u}+\cbetanorm{f}).$$
\eet

\proof We proceed as in the proof of Theorem \ref{main}, but now
in \eqref{FQE3} we write $F(x)=J_1(x)+\cdots +J_5(x)$, where
$$J_5(x)=P(x)v(x).$$
We have, using Proposition \ref{HP1} and Lemma \ref{HL1},
\begin{align*}
\cbetanorm{J_5}&\leq c_2\cbetanorm{P}\cbetanorm{\vp}\cbetanorm{u}\\
&\leq c_3 (\linorm{u}+\eps\cabnorm{u}).
\end{align*}
Other than this additional term, the rest of the proof goes through
as before.
\qed

\subsection{First order terms}\label{SS3}

If $\al>1$, we can add a first order term to $\sL$. (We can also keep
the zero order term as in Theorem \ref{FRRT1}, but we omit this
in the following discussion for simplicity.) 

\bet\label{FRRT2}
Suppose $\al>1$. For $i=1, \ldots, d$, let $Q_i$ be functions such that $\cbetanorm{Q_i}<\infty$. Let
$$ \sL'' u(x)=\sL u(x)+ \sum_{i=1}^d Q_i(x) D_iu(x),$$
where $\sL$ is defined by \eqref{I-def1} or \eqref{I-def2} and satisfies
the assumptions of Theorem \ref{main}.
Then there exists $c_1$ (which depends on $\sum_{i=1}^d \cbetanorm{Q_i}$) such that
if $\sL'' u(x)=f(x)$ and $\cabnorm{u}<\infty$, then
$$\cabnorm{u}\leq c_1(\linorm{u}+\cbetanorm{f}).$$
\eet

\proof  As in the proof of Theorem \ref{FRRT1} we have an additional
term in the definition of $F$, but this time the term is
$$J'_5(x)=\sum_{i=1}^d Q_i(x) D_iv(x).$$
We have
\begin{align*}
\cbetanorm{Q_iD_iv}&\leq c_2\cbetanorm{Q_i}(\cbetanorm{\vp D_iu}
+\cbetanorm{u D_i\vp})\\
&\leq c_3(\cbetanorm{\vp}\cbetanorm{D_iu}+\cbetanorm{u}\cbetanorm{D_i\vp})\\
&\leq c_4\linorm{u}+\eps\cabnorm{u},
\end{align*}
using Lemma \ref{HL1} and Proposition \ref{HP1}. With $J'_5$
handled in this fashion, we proceed as before.
\qed

\subsection{Higher order smoothness}\label{SS4}

One would expect that if $f$  and $A(\cdot, h)$ have additional smoothness,
then  the solution $u$ to $\sL u=f$ should have additional smoothness.
This is indeed the case. One way to show this is to extend the
estimates previously proved to $C^\beta$ and $C^\apb$ when $\beta>1$. Here is an alternate
way. We do the case $\beta\in (1,2)$ for concreteness, but the case when
$\beta\in (m,m+1)$ for some $m$ is similar. When we write $D_iA(x,h)$,
we mean  the $i^{th}$ partial derivative in the variable $x$.

\bet\label{FRRT3} Suppose $\beta\in (1,2)$ and there exists $c_1$ such that
for each $i=1, \ldots, d$, 
$$\sup_{x} \sup_h |D_iA(x+k,h)-D_iA(x,h)|\leq c_1|k|^{\beta-1}.$$
Then there exists $c_2$ such that if $f\in C^\beta$ and $u\in C^\apb$ with
$\sL u=f$, we
have
$$\cabnorm{u}\leq c_1(\linorm{u}+\cbetanorm{f}).$$
\eet

\proof We sketch the proof, and we restrict our attention to $\al<1$
for simplicity. Differentiating $\sL u=f$ yields
$$\sL (D_iu)(x)+\int [u(x+h)-u(x)] \frac{D_iA(x,h)}{|h|^{d+\al}}\, dh
=D_if.$$
Writing $G_i(x)$ for the second term on the left, 
$$\sL (D_iu)=D_if-G_i,$$
and by Theorem \ref{main},
$$\norm{D_iu}_{C^{\beta-1}}\leq c_3(\linorm {D_iu} + \norm{D_if}_{C^{\beta-1}}
+\norm{G_i}_{C^{\beta-1}}).$$
Note $\norm{D_if}_{C^{\beta-1}}\leq c_4\cbetanorm{f}$ and 
$\linorm{D_iu}\leq c_5\linorm{u}+\eps \cabnorm{u}$. Also
$\cbetanorm{u}\leq c_6\sum_{i=1}^d \norm{D_iu}_{C^{\beta-1}}$.
So the key step is to prove that
\bee\label{FRRE4}
\norm{G_i}_{C^{\beta-1}}\leq c_7\linorm{u}+\eps \cabnorm{u}.
\eee
By arguments similar to the derivation of the estimates for $J_4$ in
the proof of Theorem \ref{main} but somewhat simpler,
$$\norm{G_i}_{C^{\beta-1}}\leq c_8\norm{u}_{C^{\apb-1}}.$$
By Proposition \ref{HP1}, the right hand side
is bounded by the right hand side of \eqref{FRRE4}.
\qed

\subsection{Sharpness}\label{SS5}

Our results are sharp in several respects. For example, one might
ask if the solution $u$ to $\sL u=f$ can be taken to be in $C^{\apb+\delta}$
for some $\delta>0$ when $f\in C^\beta$. The answer is no in general. Let $\sL=\sL_0$,
where $\sL_0$ is defined by \eqref{Dpt-E10}. Let $f$ be a $C^\beta$
function that is not in $C^{\beta+\delta}$ for any $\delta$. If the
solution to $\sL u=f$ satisfied
$$\norm{u}_{C^{\apb+\delta}}\leq c_1(\linorm{u}+\cbetanorm{f}),$$
then by Corollary \ref{fd-C45}, $f= \sL_0 u$ would be in $C^{\beta+\delta}$, a contradiction.

Another question is whether one can still obtain our main estimate
\eqref{mainest} if $A(x,h)$ only satisfies 
\bee\label{FRR-E41}
\sup_x \sup_h |A(x+k,h)-A(x,h)|\leq c_1|k|^{\beta-\delta}, \qq k\in \R^d,
\eee
for some $\delta>0$. Again the answer is no in general. Let $f$ be
a function that is in $C^\beta$ but not in any $C^{\beta+\zeta}$ for $\zeta>0$. 
Let
$w$ be a function that is in $C^{\beta-\delta}$ for some $\delta\in (0,\beta)$
but not in $C^{\beta-\delta+\zeta}$ for any $\zeta>0$. Suppose also
that $w$ is bounded below by a positive constant. Let $\sL_0$
be defined as in \eqref{Dpt-E10}, and define $A(x,h)=w(x)A_0(h)$. Then 
$\sL u(x)=w(x)\sL_0 u(x)$, and $A(x,h)$ satisfies 
\eqref{FRR-E41}. Consider the solution to $\sL u(x)=f(x)$. We have
$\sL_0 u(x)=f(x)/w(x)$. If $u$ were in $C^\apb$, then 
$f(x)/w(x)=\sL_0 u(x)$ would be in $C^\beta$, a contradiction.  

\subsection{The $\linorm{u}$ term}\label{SS6}

Our main estimate \eqref{mainest} has a $\linorm{u}$ on the right hand side. 
When can one dispense with this term? 
First we give a condition where one can do so.

Suppose one considers $\sL' u(x)=f(x)$, where $\sL'$ is defined in
Theorem \ref{FRRT1} and moreover for some $\lam>0$,  $P(x)\leq -\lam$ 
for all $x$.
If $X_t$ is the strong Markov process associated to $\sL$
(that is, the infinitesimal generator of $X$ is $\sL$, for example), 
the solution to $\sL' u(x)$ is given in probabilistic terms by
$$u(x)=-\E^x \int_0^\infty e^{\int_0^s P(X_r)\, dr} f(X_s)\, ds.$$
Under the condition that $P(x)\leq -\lam$, then
$$|u(x)|\leq \E^x \int_0^\infty e^{-\lam s} |f(X_s)|\, ds
\leq \frac{1}{\lam}\linorm{f}.$$ In this case, we have the bound 
$$\linorm{u} \leq \linorm{f}/\lam\leq \cbetanorm{f}/\lam.$$

On the other hand, if there is no zero order term, there is no reason to 
expect
that 
a bound of the form
\bee\label{FR-E81}
\linorm{u}\leq c_1\cbetanorm{f}
\eee should hold when $\sL u=f$. This bound trivially fails to hold because 
$u$ plus a constant is still a solution to the equation.  

Even when we restrict ourselves to solutions that vanish at infinity,
\eqref{FR-E81} cannot hold. To see this, let $A(x,h)$ be identically
1, so that $\sL$ is the infinitesimal generator of a symmetric stable
process, let $\ol \vp$ be defined as in the beginning of Section \ref{S:F},
and let $f_r(x)=\ol \vp(x/r)$. Then $\linorm{f_r}=1$ for all $r$,
while $[f_r]_{C^\beta}\to 0$ as $r\to \infty$ for each $\beta\in (0,1)$.
On the other hand, if $u_r$ is the solution to $\sL u=f_r$,
a scaling argument shows that $|u_r(0)|=c_1r^\al\to \infty$ as $r\to \infty$.

\subsection{Future research}\label{SS10}

We mention some directions for future research. 

\begin{enumerate}
\item
{\it Interior estimates for the Dirichlet problem.}
Can one give interior estimates
for the regularity of harmonic functions (the Dirichlet problem) and
the regularity of potentials (the analog of Poisson's equation) in bounded
domains?
\item 
{\it Boundary estimates.} To obtain a satisfactory theory, one
would like estimates on harmonic functions and potentials 
in bounded domains that are valid  up to
the boundary. 
\item {\it Symmetric processes.} Suppose instead of $\sL$ one works instead
with the Dirichlet form 
$$\sE(f,g)=\int_{\R^d}\int_{\R^d} (f(y)-f(x))(g(y)-g(x))
\frac{B(x,y)}{|x-y|^{d+\al}}.$$
The generator associated to $\sE$ is the analog  of an elliptic operator in divergence form. The Harnack
inequality and H\"older regularity for harmonic functions are known
in this setting under the assumption that $B(x,y)$ is
symmetric and bounded above and below by positive constants; 
see \cite{CK03}. However if one adds some continuity conditions to $B$,
one would expect the corresponding potentials and harmonic functions 
to have additional smoothness. 
\item {\it The parabolic case.} One could 
look at the fundamental solution or heat kernel $p(t,x,y)$, which
is equivalent to looking at the transition densities of the associated
process. One would expect that if the $A(x,h)$ (and the $B(x,h)$) have some smoothness, 
say, H\"older continuous of order $\beta$,
and are bounded above and below by positive constants,
then the $p(t,x,y)$ are not only H\"older continuous in $x$ and $y$,
but will be $C^{\al+\beta}$
in each coordinate. (In the symmetric case H\"older continuity is known, but
of a smaller order.)  This question could be asked about the transition
densities in the whole space $\R^d$ and also in bounded domains.
\item {\it Variable order.} Consider operators $\sL$ of the form
\bee\label{FRR-E51}
\sL f(x)=\int [f(x+h)-f(x)-1_{(|h|\leq 1)} \grad f(x)\cdot h]
n(x,h)\, dh,
\eee
where we assume
$$\frac{c_1}{|h|^{d+\al}}\leq n(x, h)
\leq \frac{c_2}{|h|^{d+\beta}}, \qq x\in \R^d, 1\geq |h|>0,$$ 
$0<\al<\beta<2$, and some appropriate condition is imposed on 
$n(x,h)$ for $|h|\geq 1$. 
Such an operator  is  of variable order  because
if one writes it as a pseudo-differential operator, then the order
is not fixed; see \cite{jacob}. Some progress has already been made on operators
of variable order; see \cite{mosernl} and \cite{conhar} for the operators $\sL$
in \eqref{FRR-E51} and 
see \cite{nlharn} and  \cite{bkk} for non-local Dirichlet forms of variable
order. Can one give suitable assumptions on $n(x,h)$ so that harmonic
functions and potentials have additional smoothness?  
\item 
{\it Diffusions with jumps.} If we consider  operators that are the
sum of an elliptic differential operator and a non-local operator,
the same questions could be asked as for the pure jump case:
higher order derivatives, regularity up to the boundary, transition
density estimates. (The Harnack inequality was considered in \cite{foondun1}
and \cite{foondun2}.) 
\end{enumerate}

\bs

%ZZZZZZZZZZ

\medskip

\ni Richard F. Bass\\
\ni Department of Mathematics\\
\ni University of Connecticut \\
\ni Storrs, CT 06269-3009, USA\\
\ni {\it bass@math.uconn.edu}

\end{document}